%% file: wedge2.tex
\documentclass[10pt]{article}
\usepackage[english]{babel}
\usepackage{amsmath}
\usepackage{amsfonts}
\usepackage{amssymb}
\usepackage{amsthm}
\usepackage{graphics}
\usepackage{amscd}
\usepackage{epsfig}
\usepackage{color}
\input{psfig}

\newcommand{\Z}{\mathbb{Z}}
\newcommand{\R}{\mathbb{R}}
\newcommand{\N}{\mathbb{N}}

\newcommand{\C}{\mathbb{C}}

\newcommand{\E}{\mathbb{E}} 
\newcommand{\mc}{\mathcal}

\newcommand{\eps}{\varepsilon}
\newcommand{\ind}{1\!\mathrm{I}}
\renewcommand{\P}{\mathbb{P}}

\renewcommand{\H}{\mathbb{H}}

\newcommand{\lap}{\Delta\!}

\DeclareMathOperator{\End}{End}

 \DeclareMathOperator{\cotan}{cotan}
\DeclareMathOperator{\Unif}{Unif}
\DeclareMathOperator{\sgn}{sgn}
\DeclareMathOperator{\SLE}{SLE}

\title{Reflected planar Brownian motions, intertwining relations and crossing
  probabilities}
 \author{Julien Dub\'edat\footnote {Universit\'e Paris-Sud}}

\newtheorem{Prop}{Proposition}
\newtheorem{Lem}{Lemma}
\newtheorem{Cor}{Corollary}
\newtheorem{Rem}{Remark}

\begin{document}
\maketitle

\begin{abstract}
Prompted by an example arising in critical percolation, we study some 
reflected Brownian motions in symmetric planar domains and show that
 they are intertwined with one-dimensional diffusions.
In the case of a wedge, the reflected Brownian motion is
intertwined with the 3-dimensional Bessel process. This 
implies some simple hitting distributions and sheds some light on the 
formula proposed by Watts for double-crossing probabilities in 
critical percolation.
\end{abstract}

\section{Introduction and notations}

It has been pointed out in \cite {W,LSW} that the hitting distribution of a certain reflected 
Brownian motion (RBM) in an equilateral triangle was uniform. More precisely, if the Brownian motion is 
started from one corner and reflected on the two adjacent sides with the oblique reflection angle $\pi/6$
away from the normal direction toward the 
opposite side, then the Brownian motion will hit this opposite side with uniform distribution.
Their proof uses a discrete version of this result (reflected simple random walk) on a well-chosen 
triangular lattice, and an invariance principle (i.e. the reflected random walk converges to reflected 
Brownian motion). 
This result is then combined with a locality property in order to identify the law of the whole ``hull'' of this stopped reflected 
Brownian motion with that of a chordal SLE$_6$ process (that has the
same uniform hitting distribution and has also a locality property).

In \cite {S}, Smirnov proved that SLE$_6$ is the 
scaling limit of critical percolation cluster interfaces on the triangular lattice. The main step in this 
proof is the derivation of the fact that in the scaling limit, the
above-mentioned hitting distribution (if one replaces the RBM by the
percolation interface)
becomes uniform. It makes an important use of some specific features of the triangular lattice.
Note that in the case of critical percolation, the discrete hitting probabilities are not uniform, as opposed 
to those for the reflected simple random walk.

One of the motivations of the present paper is to see whether one can generalize this hitting probability 
property of reflected Brownian motion to other symmetric shapes (for instance
other wedges) and other reflection angles. It turns out that one can choose
the reflection angles in such a way that hitting probabilities of segments
orthogonal to the symmetry axis of the shape remain
uniform. Furthermore, this fact is closely related to a relation between these 
reflected Brownian motions and one-dimensional diffusions. More
precisely, the projection of the RBM on the symmetry axis is a
diffusion in its own filtration; for wedges, one recovers the 3-dimensional Bessel process. 
This provides a new example of non-trivial intertwining relations (see
for instance \cite {RP,CPY}).
Our proofs rely also on reflected simple random walks and an
invariance principle.

We shall see that the formula proposed by Watts \cite {Wa} for
double-crossing probabilities
(i.e. simultaneous top-to-bottom and left-to-right crossing of a
quadrilateral) for critical 
percolation and for which it is not clear that the SLE$_6$ approach
will confirm it, is in fact 
satisfied (in some appropriate sense) by the reflected
 Brownian motion; this will follow from
the study of time-reversed reflected Brownian motions.

This paper is organized as follows. First, we define some reflected
random walks in wedges and study their scaling limit. Then we prove
that the limiting reflected Brownian motions are intertwined with the
3-dimensional Bessel process. Making use of this fact, we determine
the time reversals of these RBMs, and discuss an analogue of Watts'
formula. Finally, we generalize some of the results to symmetric shapes (``vases''), essential replacing the 3-dimensional
Bessel process by a general one-dimensional diffusion.

\section {Invariance principle}

Let $\alpha \in (0,\frac\pi 2)$. We will consider the wedge 
$$
C_\alpha=\{z\in\C\ :\ -\alpha\leq\arg z\leq\alpha\}\subset\C
$$ 
and the rectangular lattice 
$$L_\alpha=\cos\alpha\Z+i\sin\alpha\Z.$$
We will study random walks on the graph $\Gamma_\alpha=L_\alpha\cap C_\alpha$.

\begin{figure}[htbp]
\begin{center}
\input{wedge1.pstex_t}
\end{center}
\caption{A wedge with its associated rectangular lattice} 
\end{figure}
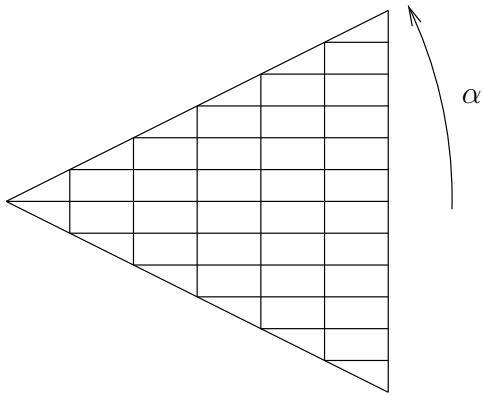

In this section, the angle $\alpha$ is fixed, so we will drop the subscript $\alpha$. Consider the following random walk on the graph $\Gamma=L\cap C$: for an inner point $x$,
the transition probability is: 
$$
\left\{\begin{array}{l}  
p=p(x,x+\cos\alpha)=p(x,x-\cos\alpha)=\frac{\sin^2\alpha}2\\
q=p(x,x+i\sin\alpha)=p(x,x-i\sin\alpha)=\frac{\cos^2\alpha}2
\end{array}\right.
$$
so that such a step has a zero mean, and its covariance matrix is a multiple of the identity.
If $x \neq 0$ is a boundary point with $\arg (x) = \alpha$, the transition probability is:
$$
\left\{\begin{array}{l}  
p(x,x+\cos\alpha)=p(x,x+\cos\alpha+i\sin\alpha)=\frac{\sin^2\alpha}2\\
p(x,x-i\sin\alpha)=p(x,x)=\frac{\cos^2\alpha}2\\
\end{array}\right.$$
Notice that in this case, 
$$\arg(\E^x(X_1-x))=2\alpha-\frac\pi 2.$$
On the other boundary line, the transition probabilities are defined symmetrically.
  
Finally for the apex, set
$$r=p(0,\cos\alpha)=p(0,\cos\alpha\pm i\sin\alpha)=\frac{1-p(0,0)}3>0.$$
The exact value of this positive probability will not matter in the continuous limit. This random walk starting from $x\in\Gamma$ will be denoted by $(X^x_n)_{n\geq 0}$. We will also call $(Y^x_t)_{t\geq 0}$ the Brownian motion in $C$ starting from $x$ and reflected on the boundary with angle $-\alpha$ (see \cite{VW} for a definition of this process; we use the same conventions for reflection angles, i.e. positive angles point towards the apex).

\begin{figure}[htbp]
\begin{center}
\input{wedge2.pstex_t}
\end{center}
\caption{Transition probabilities for an inner point and a boundary point} \end{figure}
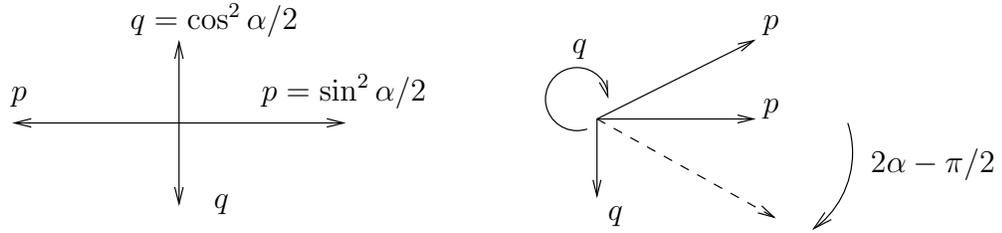

The following result is classical. We include a proof for the sake of completeness 
(and in order to point out later how it can be generalized).

\begin{Lem}
Let $x$ be a point in $C$, and $(x_n)$ a sequence in $\Gamma$ such that $|x_n-nx|\leq 1/2$ for all $n\in\N$. Then the following weak convergence of processes holds as $n$ goes to infinity: 
$$\left(\frac 1n X^{x_n}_{\lfloor
n^2t/\sin^2\alpha\cos^2\alpha\rfloor}\right)_{t\geq 0}\longrightarrow \left(Y^x_t\right)_{t\geq 0}
$$
\end{Lem}

\begin{proof} Let $(Y^{(n)}_t)$ be a continuous $C$-valued process
  that interpolates linearly the process: 
$$\left(\frac 1n X^{x_n}_{\lfloor
    n^2t/\sin^2\alpha\cos^2\alpha\rfloor}\right)$$
Choose a H\"older exponent $\beta$, $0<\beta<\frac 12$. Then it is standard to check that for any fixed $T>0$, the following estimate on $\beta$-H\"older norms holds : $$\lim_{K\rightarrow\infty}\limsup_{n\in\N}\P(||Y^{(n)}||_{[0,T],\beta}\geq
K)=0$$
One proceeds exactly as in the proof of the a.s. $\beta$-H\"older
continuity for Brownian paths (there is no particular problem due to
the boundary here). Hence we get a tight sequence in the Polish space
of $\beta$-H\"older continuous mappings from $[0,T]$ to $C$. Using
Prohorov's theorem, one gets the existence of subsequential weak
limits. Then a standard diagonal argument yields subsequential weak
limits for the whole process in the Wiener space of $C$-valued
processes.

So we have to check that there is only one possible weak limit, namely
the reflected Brownian motion $(Y_t)$. To prove this, one can use the
submartingale problem characterization of reflected BM (\cite{VW},
Theorem 2.1). More precisely, the RBM is the only $C$-valued process
starting from $x$ such that for any $f\in C^2_b(C)$ (twice
differentiable with bounded derivatives) with positive derivatives on
the boundary along the reflection direction,
 the real valued process: $$f(Y_t)-\frac 12\int_0^t\lap f(Y_s)ds$$ is
 a submartingale. So let $f$ be such a function; it is sufficient to
 check that for any $0\leq s<t$: 
$$\liminf_{n\in N}\E_{(n)}\left(f(Y_t)-f(Y_s)-\frac
    12\int_s^t\lap
f(Y_u)du\right)\geq 0$$
where $\E_{(n)}$ designates the expectation operator for the $n$-th
approximate process. Consider the stopping time $\tau=\inf_{t\geq
  0}(|Y^x_t|\geq M)$ for some large number $M$. It is enough to check
the submartingale inequalities up to time $\tau$. Notice that  
the following Taylor expansion holds: 
\begin{align*} \left(pf(x+\cos\alpha/n)+pf(x-\cos\alpha/n)+qf(x+i\sin\alpha/n)+qf(x-i\sin\alpha/n)\right)\\
=\frac
1{2n^2}\sin^2\alpha\cos^2\alpha\lap f(x)+o(\frac1{n^2}) \end{align*} and that the error term is uniform in $x\in C\cup D(0,M)$ (using for instance Hadamard's lemma). Let $u_k=k\sin^2\alpha\cos^2\alpha/n^2$,
$k\in\N$, correspond to the discrete jump times, and let $B$ be the
boundary of $C$. Then: 
\begin{align*}
\E_{(n)}^x(f(Y_{t\wedge\tau})-f(Y_{s\wedge\tau})) &= \E_{(n)}^x\left(\sum_{s\wedge\tau\leq
    u_k< t\wedge\tau}f(Y_{u_{k+1}})-f(Y_{u_k})\right) \\  
&= \E_{(n)}^x\left(\sum_{s\wedge\tau\leq
    u_k< t\wedge\tau}\frac 1{2n^2}\cos^2\alpha\sin^2\alpha\lap
  f(Y_{u_k})\right.\\
&\left.\vphantom{\sum_{u_k}}\ \ \ \ \ \ \  +\eps(Y_{u_k})\ind_{Y_{u_k}\in B}\right)+o(1)\\
&= \E_{(n)}^x\left(\int_{s\wedge\tau}^{t\wedge\tau}\frac 12\lap
  f(Y_u)du\right)+o(1)+\sum_{y\in B}\E_{(n)}^x(L_y)\eps(y)
\end{align*}
 The error term $o(1)$ is simply a Riemann sum error. Then
there is a ``boundary error term''. In the formula, $L_y$ designates
the local time at $y$ (number of discrete jumps to $y$) between times
$s\wedge\tau$ and $t\wedge\tau$, and $\eps(y)$ is defined as (if
$\arg(y)=\alpha$ for instance):
$$\eps(y)=q(f(y)-f(y+i\sin\alpha/n))+p(f(y+e^{i\alpha}/n)-f(y-\cos\alpha/n))$$
The first order term is proportional to the derivative of $f$ at $y$
along the reflection direction, hence is positive by hypothesis. One
may write: $\eps(y)+O(1/n^2)\geq 0$ for $y\in B$, the error term
$O(1/n^2)$ being uniform in $y\in C\cup D(0,M)$. So the only thing to
check now is that the time spent on the boundary is negligible, more
precisely:
$$\sum_{y\in B}\E_{(n)}^x(L_y)=o(n^2)$$
Recall that
we consider the local time between times $s$ and $t$, and before exit time $\tau$; one may drop this last condition for the sake of simplicity. Justifying the very intuitive fact that the walk spends a negligible time on a negligible part of the state space is rather tedious. The discrete intertwining relations allow explicit computations that yield the result when the starting point is the apex $0$ (see below). Then, using the Markov property, one sees it is enough to prove this estimate for zero. This concludes the proof. \end{proof}

One may carry out the same proof with different transition probabilities on the boundary, corresponding to different reflection angles; to derive the final local time estimate in the general case, one may use the local Central Limit Theorem.

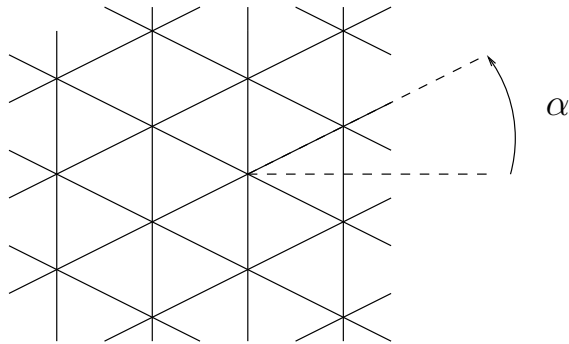
\begin{figure}[htbp]
\begin{center}
\input{wedge3.pstex_t}
\end{center}
\caption{``stretched'' triangular lattice} \end{figure}

\begin{Rem}
Similar results hold for different lattices, in particular ``stretched triangular lattices'', i.e. images of the standard triangular lattice under an orthogonal affinity. In the special case $\alpha=\pi/6$, one may use the standard triangular lattice, with the usual transition probability for inner point (i.e., $1/6$ for each neighbour). In the case $\alpha=\pi/4$, these two approaches give the same lattice, namely the square lattice. 
\end{Rem}

Note that all these reflected Brownian motions have the same parameter $\alpha=-1$ in Varadhan-Williams conventions (\cite{VW}) (what we call $\alpha$ would be denoted $\xi/2$ in \cite{VW}).

\section{Intertwining Relations}

First, let us recall the notion of intertwined Markovian semigroups
(for some background on the subject, see \cite{CPY}). Let $(P_t,t\geq
0)$ and $(Q_t,t\geq 0)$ be two Markovian semigroups with respective
state spaces $(S_P,\mc S_P)$ and $(S_Q,\mc S_Q)$, and $\Lambda$ a
Markov transition kernel from $(S_Q,\mc S_Q)$ to $(S_P,\mc S_P)$. The
two semigroups are said to be intertwined by $\Lambda$ if the
following identity of Markov transition kernels from $S_Q$ to $S_P$
holds : $$\Lambda P_t=Q_t\Lambda$$ for all $t\geq 0$ ($t$ is either a
discrete or a continuous time parameter).

Suppose that there exists a measurable function $\phi: S_P\rightarrow
S_Q$ such that $\Lambda(y,\phi^{-1}(y))=1$ for all  $y\in S_Q$. Then, as shown in \cite{RP}, if $(Z_t)$ is a Markov process with semigroup $(P_t)$ and initial law $\Lambda(y_0,.)$ for some $y_0\in S_Q$, and $Y_t=\phi(Z_t)$, then $(Y_t)$ is a Markov process with semigroup $(Q_t)$ and starting state $y_0$. Moreover, for any fixed time $T$, the following filtering formula holds: 
$$
\P(Z_T\in .|Y_t,t\geq 0)=\Lambda(Y_T,.).
$$ 
Under appropriate regularity conditions (which will be satisfied in all examples that we shall consider here), this formula also holds for almost sure finite stopping times $T$ (w.r. to the filtration of $Y$). In particular, for $T=T_y$, the first hitting time of $y$ by $Y$, if $T_y<\infty$ a.s., then the random variable $Z_{T_y}$ has distribution $\Lambda(y,.)$ and is independent of $(Y_t)$. 

Consider as before the random walk $(X_n)$ on the graph $\Gamma$. Let $\Lambda$ be the following Markov transition kernel from $\N$ to $\Lambda$ (we omit to mention that these discrete spaces are equipped with  discrete
$\sigma$-algebras):
$$\Lambda(x,.)=\frac 1{2x+1}\sum_{k=-x}^x\delta_{x\cos\alpha+ki\sin\alpha}$$
Let $(P_n)$ be the semigroup of $(X_n)$, and let $(Q_n)$ be the
semigroup of the random walk on $\N$ with transition probability:
 $$\left\{\begin{array}{l} q(i,i-1)=\frac{\sin^2\alpha}2\frac
     {2i-1}{2i+1}\\ q(i,i)=\cos^2\alpha\\ q(i,i+1)=
     \frac{\sin^2\alpha}2\frac {2i+3}{2i+1} \end{array}\right.$$
 for $i>0$ and $q(0,0)=1-3r,\ q(0,1)=3r$. Note that $i\mapsto \frac 1{2i+1}$ is harmonic for this transition kernel, except at 0. Then the intertwining relation
holds:
$$\Lambda P_n=Q_n\Lambda$$
for all $n\geq 0$.

Let $(Y_t)$ be the reflected Brownian motion in the wedge $C$ with
reflection angle $-\alpha$. Let $(P_t)_{t\geq 0}$ be its
semigroup. Let $(Q_t)_{t\geq 0}$ be the generator of the 3-dimensional
Bessel process (with values in $\R_+$), and $\Lambda$ be the Markov
transition kernel from $\R_+$ to $C$ such that $\Lambda(x,.)$ is the
uniform distribution on the segment $[x-i\tan\alpha,x+i\tan\alpha]$.

In the scaling limit, the previous intertwining leads to:

\begin{Prop} 
Then $P$ and $Q$ are intertwined by $\Lambda$. 
In particular, for $x\geq 0$, $Z_t=\Re Y^{\Lambda(x,.)}_t$ 
is a 3-dimensional Bessel process starting from $x$. 
\end{Prop}

\begin{proof}
It is obvious that the intertwining relation holds in the limit. The only thing to check is that the intertwined random walk on $\N$ converges to a 3-dimensional Bessel process. One may proceed as in the previous proof, after computing: \begin{eqnarray*} 
\left(\cos^2\alpha
  f(x)+\frac{\sin^2\alpha}2\left(\frac{2nx/\cos\alpha+3}{2nx/\cos\alpha+1}f(x+\cos\alpha/n)\right.\right.\\
+\left.\left.\frac{2nx/\cos\alpha-1}{2nx/\cos\alpha+1}f(x-\cos\alpha/n)\right)\right)-f(x)
 \\
 =
\frac{\sin^2\alpha\cos^2\alpha}{n^2}\left(\frac 1xf'(x)+\frac 12f''(x)\right)+o(1/n^2) 
\end{eqnarray*}
Here one recognizes the generator of the 3-dimensional Bessel process
(the factor $\frac{\sin^2\alpha\cos^2\alpha}{n^2}$ corresponds to the
time scaling). Such discrete approximations of the 3-dimensional
Bessel process are classical.
\end{proof}

The Markov transition kernel $\Lambda$ acts on bounded Borel functions
in the following way: if $f$ is a bounded Borel function on $C$, then
$\Lambda f$ is a bounded Borel function on $\R_+$ such that:
$$(\Lambda f)(x)=\E(f(x+iU\tan\alpha))$$ where $U$ is a uniform random variable on $[-1,1]$. One may remark that these intertwining relations fit in the ``filtering type framework'' described in \cite{CPY}.

As Jim Pitman pointed out to us, there is an analogy with the
situation in  the $(2M-X)$ theorem: let $B$ be a standard (real)
Brownian motion starting from 0, and let $M_t$ be its supremum up to
time $t$. Then the process $(2M_t-X_t)_{t\geq 0}$ is a 3-dimensional Bessel
process (in its own filtration) (see \cite{RP,RY}). A proof of this fact involves
an intertwining relation. More precisely, note $X_t=M_t-B_t$,
$Y_t=2M_t-B_t$. Then, if $({\mc G}_t)_{t\geq 0}$ designates the natural
filtration of $(Y_t)_{t\geq 0}$, the following relation holds for every Borel
function $f: \R_+\rightarrow\R_+$:
$$\E(f(X_t)|{\mc G}_t)=\int_0^1 f(xY_t)dx$$
The Markov kernel in this last situation: $\Lambda f(y)=\int_0^1 f(xy)dy$
is very similar to the one we described earlier.

The Brownian motion in a plane strip $\{z:\ 0\leq \Im z\leq 1\}$ with
normal reflection on the boundary may be seen as a degenerate limiting
case when $\alpha\searrow 0$. Indeed, the reflection direction makes
an angle $\pi/2-\alpha$ with the boundary, and the two boundary
half-lines make an angle $2\alpha$. This particular RBM may be
represented as $Z=X+if(Y)$, where $X$ and $Y$ are independent standard
(one-dimensional) BMs, and $f$ is the ``seesaw'' function: for all
$k\in\Z$, 
\begin{align*}
f(x)=x-2k&&\text{if\ \ \ \ }&2k\le x\le 2k+1\\
f(x)=2k+2-x&&\text{if\ \ \ \ }&2k+1\le x\le 2k+2
\end{align*}
In this case, the intertwining relation
is easily proved. Indeed, let $\phi$ be any bounded Borel function on
the strip and $t$ be a fixed positive time. We have to check that:
$$\E_x(\phi(X_t+iU))=\E(\E_{x+iU}(\phi(Z_t)))$$
where $U$ is an independent uniform random variable. Since the real
and the imaginary part of $Z$ are independent, we have only to check
that, for any bounded Borel function $\varphi$ on $[0,1]$:
$$\E(\varphi(U))=\E(\E_U(\varphi(f(Y_t))))=\E_0(\varphi(f(U+Y_t)))$$ 
which is readily seen (if $V$ is any symmetric r.v. independent from $U$,
$f(U+V)$ is a uniform r.v.).

\begin{Rem} In the general case, the intertwining relation does not
  seem obvious directly in the
  continuous setting, hence the use of discrete approximations. 
  Also, the real part of the RBM started from a
  fixed, inner point (not the apex) does not appear to have similar
  properties; in this case, the relation between the filtration of the
  RBM and the filtration of its real part seems very intricate, and
  the real part is generically not Markovian in its own filtration.
\end{Rem}

\section{Time reversal}

In this section, we consider the time reversals of the reflected Brownian motions studied in the previous sections. The intertwining relation enables an explicit computation of discrete Green functions, hence the determination of the time reversal of discrete random walks; then one takes the continuous limit.

\begin{Prop}
The time reversal of the RBM with reflection angle $-\alpha$ starting from $0$ and stopped when hitting $\{z:\Re z\geq x\}$ is the RBM with reflection angle $+\alpha$ starting from the uniform distribution on the segment $[x-i\tan\alpha;x+i\tan\alpha]$, conditioned not to hit $\{z:\Re z\geq x\}$ again, and killed at $0$. This last process is intertwined with a 3-dimensional Bessel process. 

\end{Prop}

Loosely speaking, the time reversal of a reflected Brownian motion in $C_\alpha$ with reflection angle $-\alpha$ starting from $0$ is the reflected Brownian motion reflected at angle $+\alpha$, killed at $0$, starting from ``the uniform distribution at infinity''. 

\begin{proof}
Recall that we considered a random walk $X$ on $\Gamma$ intertwined
with a random walk $Y$ on $\N$. Let $G$ be the Green function for $X$
killed when it hits $\{z:\Re z\geq N\cos\alpha\}$ for some large
$N$. Let $\tilde G$ be the Green function for $Y$ killed at level
$N$. Then it is obvious from the intertwining of $X$ and $Y$ that:
$$G(0,x\cos\alpha+iy\sin\alpha)=\frac1{2x+1}\tilde G(0,x)$$ where
$0\leq x\leq N$, $-x\leq y\leq x$, $x,y\in\N$. As previously
mentioned, the function $i\mapsto \frac 1{2i+1}$ is harmonic for the
Markov kernel of $Y$ (except in $0$). From this, it is easy to
compute: $$\tilde
G(0,y)=\frac{(2y+1)(1-\frac{2y+1}{2N+1})}{\cos^2\alpha}$$
Then Nagasawa's formula (see e.g. \cite{RW}, III.42) yields the Markov
kernels of the time reversals $\hat X$ and $\hat Y$ of $X$ and
$Y$. For instance, we record the transition probability for $\hat X$,
if $x\in\Gamma$ is an inner point: $$\left\{\begin{array}{l} \hat
    p(x,x+i\sin\alpha)=\hat p(x,x-i\sin\alpha)=\frac{\cos^2\alpha}2\\
\hat p(x,x-\cos\alpha)=\frac{\sin^2\alpha}2\frac{N-x+1}{N-x}\\
\hat p(x,x+\cos\alpha)=\frac{\sin^2\alpha}2\frac{N-x-1}{N-x}
\end{array}\right.
$$
Notice that, as $N$ goes to infinity, one gets the original transition
probabilities. On the boundary, the reflection angle is reversed, which
is not surprising.
 At this point, one takes the continuous limit as in section 1 (obviously the time reversal operation is compatible with the continuous limit). The various claims follow easily. 
\end{proof}

The conformal equivalence $z\mapsto z^\beta$ maps $C_\alpha$ onto $C_{\beta\alpha}$. After an appropriate time change, this yields a general result for the time reversal of a reflected Brownian motion in an infinite wedge.

In the set-up of the proposition, one may notice that up to its first
hitting of the boundary, the time-reversed RBM has the law of a
Brownian excursion in the half-plane $\{z:\ \Re z\leq x\}$ (for
background on Brownian excursions, see \cite{LSW,V}). One may use this
to compute some probabilities of (indirect) interest in critical
percolation. For simplicity, we will consider a RBM $Z$ in the cone
$\Delta=(1/2+i\sqrt 3/2)-iC_{\pi/6}$ starting from the apex and
stopped on hitting the real line at time $\tau$ (we are looking at a
RBM in an equilateral triangle). We have seen that $Z_\tau$ is
uniformly distributed on $[0,1]$. Let $g$ be the last time spent by
$Z$ on the boundary before $\tau$. We are interested in the joint law
$(Z_g,Z_\tau)$, or rather in which side was last visited by $Z$
conditionally on the exit point $Z_\tau$; so we will consider the event: 
$${\mc R}=\{\Re Z_g\geq 1/2\}$$
Since the time reversal of $Z$ is a Brownian excursion until it hits the boundary, we have to
compute the ``harmonic measure'' for the Brownian excursion.   

Let $Y$ be a Brownian excursion starting in $\H=\{z: \Im z\geq 0\}$,
and let $T$
 be the first time it hits $\partial\Delta$. Let $\eps>0$, $M>0$, and
 $B$
 a complex Brownian motion exiting the strip $\{z:\ 0<\Im z<M\}$ at
 time
 $T_M$. As $M$ goes to infinity, the Brownian motion $B$ conditioned
 on exiting
 the strip by the top converges to the Brownian excursion.
 Hence, making use of the Markov property of $B$: \begin{align*} \P_{x+i\eps}(Y_T\in dy)&=\lim_{M\rightarrow\infty}\P_{x+i\eps}(B_T\in
dy|\Im B_{T_M}=M)\\
&=\P_{x+i\eps}(B_T\in dy)\lim_{M\rightarrow\infty}\frac{\P_y(\Im
  B_{T_M}=M)}{\P_{x+i\eps}(\Im B_{T_M}=M)}\\
&=\frac{\Im y}\eps\P_{x+i\eps}(B_T\in dy)
\end{align*}
It is well known that harmonic measure is conformally invariant
(see e.g.\\ \cite{RY}),
 and that the harmonic measure on the real line seen from $a+ib$ is
 given
 by a Cauchy law (see e.g. \cite{RW}) with density:
 $$\frac {bdx}{\pi(b^2+(x-a)^2)}$$
 Moreover, according to the Schwarz-Christoffel formula (see
 \cite{Ahl}), the holomorphic map: 
$$F(z)=\int_0^zu^{-2/3}(1-u)^{-2/3}du/ B(1/3,1/3)$$
 is the conformal equivalence between the upper half-plane $\H$ and
 the equilateral triangle
 $\Delta\cap\H$ that maps $(0,1,\infty)$ to $(0,1,1/2+i\sqrt 3/2)$.
 If $a\in (0,1)$ is such that $F(a)=x$, we see that:
$$P_x(Y_T\in dy)=\Im y\frac{dF^{-1}(y)}{\pi F'(a)(F^{-1}(y)-a)^2} $$
 Then, one may compute:
 \begin{align*} 
\P({\mc R}| Z_\tau\in dx) &=F'(a)^{-1}\int_1^\infty
\frac{dt}{\pi(t-a)^2}\int_1^t
\frac{(u(u-1))^{-2/3}}{B(1/3,1/3)}\frac{\sqrt 3}2du\\
&=F'(a)^{-1}\frac{\sqrt 3}{2\pi B(1/3,1/3)}\int_1^\infty\frac{du}{(u(u-1))^{2/3}(u-a)}\\
 &=\frac{\pi\sqrt 3}{3\Gamma(2/3)^3}(a(1-a))^{2/3}
\ _2F_1(1,4/3;5/3;a)\\
&=\frac{1}{B(2/3,2/3)}\int_0^a\frac{dt}{(t(1-t))^{1/3}}
 \end{align*}
 where $_2F_1(1,4/3;5/3;.)$ designates a generalized hypergeometric
 function (see eg \cite{B} or \cite{AS}, especially formulas 15.3.1
 and 15.2.5). 

\section{Relation with Watts' formula}

In this section, we recall Watts' formula (\cite{Wa}) and explain how it may be
translated in the SLE language. This relation was suggested to us by
Wendelin Werner.
We shall see that it is closely related with the formula we
derived above; in fact, this is one of the 
initial motivations for the present paper.

Recall that Cardy's formula \cite {Ca} gives the asymptotic behaviour
of the probability of an open crossing between two sides of a 
topological rectangle in the limit when the mesh of the lattice goes
to zero (for critical 
percolation). This is equivalent to the fact that the hitting
distribution of the exploration 
process in an equilateral triangle is uniform. This was shown by
Smirnov \cite {S} to hold in the case 
of site percolation on the triangular lattice, and it follows that the
scaling limit of the 
whole exploration process is $\SLE_6$. Note that it is not difficult
(see \cite{LSW1}) to prove directly in the continuous setting that this
hitting distribution for $\SLE_6$ is uniform.

Trying to generalize Cardy's results using the same approach based on
conformal field theory, Watts
\cite {Wa} considered the event that there exists simultaneously an
open left-to-right crossing and an open 
top-to-bottom crossing of a topological rectangle. He proposed a
formula to describe the asymptotic
behaviour of the probability of this event when the mesh of the
lattice vanishes, that seems to fit
well with numerical simulations.
Just as in the case of Cardy's formula, the double-crossing event can
in fact be rephrased in terms of the 
exploration process. 

For simplicity, we will discuss critical percolation on the triangular
lattice (so each vertex is colored in black with probability $p=1/2$).
Let $R$ be a topological rectangles, its boundary consisting
of four disjoint arcs: $\partial_L$, $\partial_T$, $\partial_R$,
$\partial_B$ (left, top, right, bottom). Let
$C(\partial_i,\partial_j,c)$ be the event that there exists a crossing
between $\partial_i$ and $\partial_j$ with color $c$ (here
$i,j\in\{L,T,R,B\}$, and the color $c$ is black or white:
$c\in\{b,w\}$). We will also need the event $T(\partial_i,c)$ that
there exits a connected component with color $c$ linking the four
boundary arcs except maybe $\partial_i$. 
As the triangular lattice is
self-matching (see \cite{Ke}), it is classical that the two events
$C(\partial_L,\partial_R,w)$ and $C(\partial_T,\partial_B,b)$ are
complementary:
$$C(\partial_L,\partial_R,w)=\overline{C(\partial_T,\partial_B,b)}$$
Now we are interested in $C(\partial_L,\partial_R,b)\cap
C(\partial_T,\partial_B,b)$. With a little plane topology, one sees
that:
\begin{align*}
C(\partial_L,\partial_R,b)\cap
C(\partial_T,\partial_B,b)=T(\partial_T,b)\backslash
C(\partial_L,\partial_R,w)\\
T(\partial_T,b)=C(\partial_L,\partial_R,b)\backslash T(\partial_T,w)
\end{align*}
Figure $4$ illustrates these relations.

\begin{figure}[htbp]
\begin{center}
\input{wedge4.pstex_t}
\end{center}
\caption{Paths in a topological rectangle}
\end{figure}
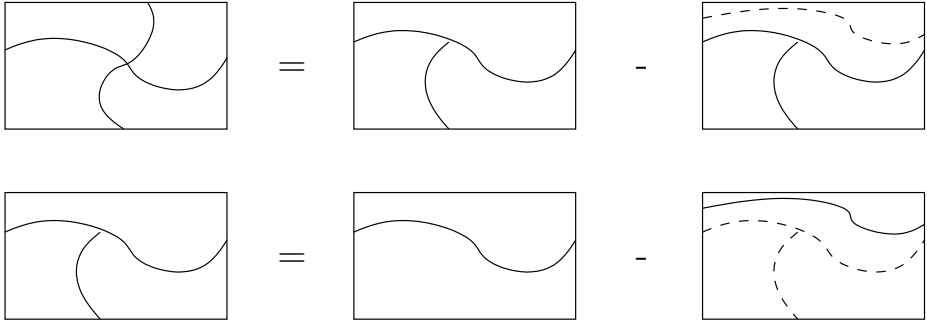

Translating this into probabilities:
$$
\P(C(\partial_L,\partial_R,b)\cap
C(\partial_T,\partial_B,b))=\P(T(\partial_T,b))-\P(T(\partial_T,b)\cap
C(\partial_L,\partial_R,w))
$$
and
$$
\P(T(\partial_T,b))=\P(C(\partial_L,\partial_R,b))-\P(C(\partial_L,\partial_R,b)\cap
T(\partial_T,w))
$$
Keeping in mind that each vertex is colored with probability $1/2$, so
that changing all the colors is measure-preserving, one gets:
$$\P(C(\partial_L,\partial_R,b)\cap
C(\partial_R,\partial_L,b))=\P(C(\partial_L,\partial_R,b))-2\P(C(\partial_L,\partial_R,w)\cap
T(\partial_T,b))$$

Consider now the situation in an equilateral triangle $ABC$. Let $X$
be the rightmost point on $BC$ that is separated from $AC$ by a black
path. Cardy's formula tells that $X$ is uniformly distributed on
$BC$ (see figure). The topological rectangle $ADXE$ delimited by the upper-half of $AB$, the
rightmost black path (solid), the leftmost white path (dashed), and
the upper-half of $AC$ is either crossed by a white path from $AD$ to
$XE$ or by a black path from $DX$ to $EA$. The exploration path goes
from $A$ to $X$ leaving white vertices on its left and black vertices on
its right. In the case where there is a white crossing between $AD$
and $XE$, it is clear that the last edge of $ABC$ visited by the
exploration process before it reaches $X$ is $AB$. Conversely, if
there is a black crossing between $DX$ and $EA$, the exploration
process last visits $AC$. 

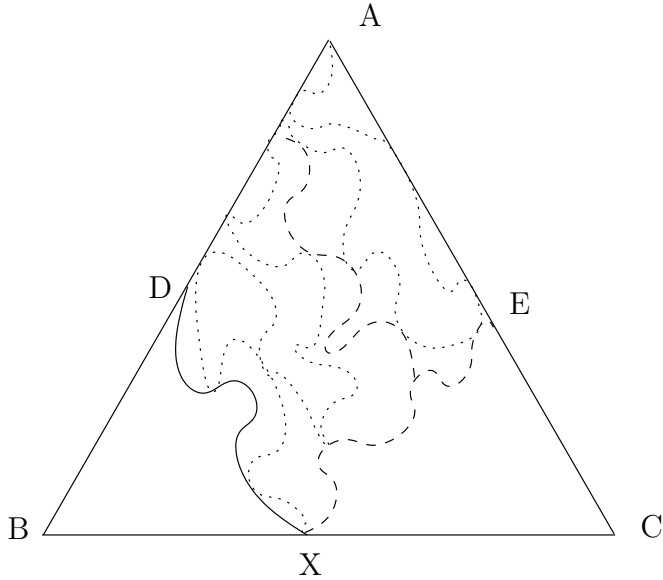
\begin{figure}[htbp]
\begin{center}
\input{wedge5.pstex_t}
\end{center}
\caption{The -stylized- exploration process (dotted)}
\end{figure}

Smirnov has proved that the exploration process converges to a chordal
$\SLE_6$ (see \cite{S}). So let $\gamma$ denote the trace of a
chordal $\SLE_6$ from $A$ to $B$ in $ABC$, $\tau$ be the first time it
hits $BC$, and $g$ the last time it hits $(AB)\cup(AC)$ before
$\tau$. We may now state the $\SLE$ interpretation of Watts' formula
in the topological rectangle $ABXC$ (the left-hand should be
understood as the scaling limit of the corresponding discrete probabilities):
$$\P(C(AB,XC,b)\cap
C(BX,CA,b))=\P(\gamma_\tau\in BX)-2\int_0^X\P(\gamma_g\in AC|\gamma_\tau=Y)d\P(\gamma_\tau=Y)$$
Recall that $\P(\gamma_\tau\in BX)=BX/BC$ (Cardy's formula). Watts'
conjectural formula is equivalent to:
$$\P(\gamma_g\in AC|\gamma_\tau=X)=\frac{1}{B(2/3,2/3)}\int_0^{BX/BC}\frac{dt}{(t(1-t))^{1/3}}$$
Let us explicit the relation between chordal $\SLE_6$ in $ABC$ and the
$RBM$ we studied in the previous section. As mentioned in the
introduction, the hull of chordal $\SLE_6$ in $ABC$ from
$A$ to $B$ stopped at time $\tau$ and the hull of the RBM in $ABC$
stopped at its first hitting of $BC$ are identical in law. Though, the
processes producing these hulls are starkly different (see
\cite{LSW}). Indeed, the $\SLE$ trace is non self-traversing, while
the RBM is likely to go back through its past hull. Let $Z$ be such a RBM,
$\tau$ the first time it hits $BC$, $g$ the last time before $\tau$ it
visits $(AB)\cap (AC)$. We have seen that:
$$\P(Z_g\in AC|Z_\tau=X)=\frac{1}{B(2/3,2/3)}\int_0^{BX/BC}\frac{dt}{(t(1-t))^{1/3}}$$
and this is {\it exactly} the formula proposed by Watts for the
scaling limit of critical percolation clusters, i.e. $\SLE_6$. Note also
that the lowest points of the hulls on $AB$, say, have the same law
since they are hull measurables. Since $\gamma$ is non
self-traversing, $\gamma_g$ is either the lowest point on $AB$ or the
lowest point on $AC$; but a.s. $Z_g$ is not the lowest point on $AB$
or $AC$. 


Given these differences between the  underlying curve laws, replacing the $\SLE_6$ by the reflected Brownian motion in the ``last visited edge'' probability 
problem, the answer should a priori be different. 
Hence, there are two possibilities:
\begin {itemize}
\item 
For some reason, the $\SLE_6$ and RBM ``last visited edge'' probabilities coincide, and Watts' prediction 
holds.
\item
Watts' prediction does not hold.
\end {itemize}

\section{Vases}

In this section, we generalize the previous properties of RBM, replacing
wedges by more general domains with a symmetry axis. More precisely,
let $h:\R_+\rightarrow\R_+$ be a positive, differentiable
function. Suppose $h(0)=0$, and $h(x)>0$ for $x>0$. Then we will
consider the ``vase'' $C_h$: 
$$C_h=\{z\in\C\ :\ \Re z\geq 0, |\Im
z|\leq h(\Re(z))\}$$ 
The shape function $h$ is fixed in this section, so we may omit the subscript.
As in the previous section, we define a tailor-made discrete model to prove an intertwining relation for its continuous limit. The situation here is slightly more complicated, so we will use a continuous time, discrete state space Markov chain. 
For the sake of simplicity, suppose the function $h$ is strictly increasing. Let $N$ be some large, fixed integer, and let $x_k=h^{-1}(k/N)$. We will consider a graph $\Gamma$ (depending on $h$, $N$) with vertices $(x_k+iy/N)_{k\in\N,y\in\Z,|y|\leq k}$. For $k\in\N$, define $\alpha_k=\arctan(1/(x_{k+1}-x_k))$. Let us now define the jump rates (non diagonal elements of the $Q$-matrix). If $z=x_k+iy/N$ is an inner point of $\Gamma$ (i.e $|y|<k$), the jump rates are : $$\left\{\begin{array}{l} 
Q(z,z+i)=Q(z,z-i)=\frac 12\\
Q(z,x_{k+1}+iy/N)=[\cotan\alpha_k(\cotan\alpha_k+\cotan\alpha_{k-1})]^{-1}\\

Q(z,x_{k-1}+iy/N)=[\cotan\alpha_{k-1}(\cotan\alpha_k+\cotan\alpha_{k-1})]^{-1}\\
\end{array}\right.$$
These jump rates are chosen in order to get a zero mean, isotropic walk. For a boundary point $z=x_k+ik/N$: $$\left\{\begin{array}{ll} 
Q(z,z-i/N)&=1/2\\
Q(z,x_{k+1}+ik/N)&=Q(z,x_{k+1}+i(k+1)/N)\\
&=[\cotan\alpha_k(\cotan\alpha_k+\cotan\alpha_{k-1})]^{-1}
\end{array}\right.$$
For the apex, set $Q(0,x_1)=Q(0,x_1\pm 1/N)>0$. Here the jump rates
are chosen so that intertwining relation holds. More precisely, let $\Lambda$ be the Markov transition kernel from $(x_k)_{k\in\N}$ to $\Gamma$ defined as: $$\Lambda(x_k,.)=\frac 1{2k+1}\sum_{y=-k}^k\delta_{x_k+iy/N}$$
Then the intertwining relation is valid for $Q$-matrices. More
precisely, let $Q$ be the $Q$-matrix of the Markov process just
defined, and $\tilde Q$ be the $Q$-matrix of the Markov process on
$\{x_k\}_{k\in\N}\subset\R_+$ defined by the jump rates:
 $$\left\{\begin{array}{l} \tilde{Q}(x_k,x_{k+1})=\frac{2k+3}{2k+1}[\cotan\alpha_k(\cotan\alpha_k+\cotan\alpha_{k-1})]^{-1}\\
\tilde{Q}(x_k,x_{k-1})=\frac{2k-1}{2k+1}[\cotan\alpha_{k-1}(\cotan\alpha_k+\cotan\alpha_{k-1})]^{-1}
\end{array}\right.$$
for $k>0$, and $\tilde Q(0,x_1)=3Q(0,x_1)$. It is immediate to check that: $$\Lambda Q=\tilde Q\Lambda$$ 
Then it is sufficient to exponentiate this intertwining relation (in
appropriate Banach spaces, say $\End(l_1(\Gamma))$ and
$\End(l_1(\{x_k\}_{k\in N}))$\ ) to get the intertwining relation for the associated semigroups $(P_t)$, $(\tilde P_t)$: $$\Lambda P_t=\tilde P_t\Lambda$$ 
for all positive time $t$. We now examine the scaling limit of the second process. Let $f$ be some arbitrary function in $C^2_b(\R_+)$. Then, as $N$ goes to infinity, if $k=k_N$ is such that $x_{k_N}\rightarrow x$ for some $x\in\R_+$: \begin{align*}
Qf(x_k)
&=\frac 1{(2k+1)(\cotan\alpha_k+\cotan\alpha_{k-1})}\left(\frac{2k+3}{\cotan\alpha_k}
    (f(x_{k+1})-f(x_k))+\right.\\
&\ \ \ \ \ \ \ \ \ \ \ \ \ \left.\frac{2k-1}{\cotan\alpha_{k-1}}(f(x_{k-1})-f(x_k))\right)\\
&=\frac 1{N^2}\left(\frac{h'(x)}{h(x)}f'(x)+\frac 12
  f''(x)\right)+o\left(\frac 1{N^2}\right)
\end{align*}
After a time rescaling, these discrete generators ``converge'' to the
diffusion generator: 
$$\frac{h'}{h}\frac\partial{\partial x}+\frac
12\frac{\partial^2}{\partial x^2}$$
 One checks immediately that a scale
  function for this diffusion on $\R_+$ is: $$\phi(x)=\int_1^x\frac 1{h^2(u)}du$$ For a wedge, $h(x)=x\tan\alpha$, then $x\mapsto x^{-1}$ is a scale function for the associated diffusion, i.e. the 3-dimensional Bessel process.

We now discuss the proof for the convergence of the processes defined
above. It is similar to the previous proof, with added
technicalities. Firstly, one has to adapt the submartingale-problem
characterization for vases. One possibility is to map conformally
$C_h$ to a wedge, say $C_{\pi/4}$ (for the Riemann Mapping Theorem,
see \cite{Ahl}); after a time change (see \cite{RY} for a discussion
of the image of a complex Brownian motion under a conformal map), one
may use a variant of the Varadhan-Williams result. Here the reflection
angle may vary along the boundary, but stay negative. Then one argues
as above: subsequential scaling limits exist (tightness); they satisfy
a submartingale-problem -a martingale problem for the one-dimensional
diffusion-, as seen by dominated convergence after applying
Dynkin's formula; hence the limit is uniquely determined, and the
processes converge weakly. Reflections are taken care of as before,
via an occupation time estimate. The intertwining is obviously
preserved in the continuous limit. Hence one may state:

\begin{Prop}
The reflected Brownian motion in the vase $C_h$ with constant reflection index $-1$ is intertwined with the diffusion on $\R_+$ corresponding to the generator  $\frac{h'}{h}\frac\partial{\partial x}+\frac 12\frac{\partial^2}{\partial x^2}$ by the Markov transition kernel from $\R_+$ to $C_h$: $$\Lambda(x,.)=\Unif([x-ih(x);x+ih(x)])$$
\end{Prop}

Let us stress that the reflection angle varies along the boundary; the
reflection at $x+ih(x)$ points in the $-i\exp(2i\arctan h'(x))$
direction.
Starting from an arbitrary diffusion on $\R_+$, one may consider an
increasing scale function $\phi$. Setting $h(x)= 1/\sqrt{\phi'(x)}$,
one constructs an intertwined reflected Brownian motion. Of course,
one can also consider shape functions defined and positive on $\R$
(``funnels'' rather than vases), or on a segment $[a,b]$, with
$h(x)=0$ if and only if $x\in\{a,b\}$ (closed vases). 
We will now consider an example of the previous construction.

\begin{Cor}
Consider the shape function $h(x)=x^\beta$, with $\beta>0$. Then the
real part of the RBM in the vase $C_h$ with constant reflexion index
$-1$ is a Bessel process of dimension $2\beta+1$.
\end{Cor}

When the shape function $h$ is not monotone, one has to be a bit more cautious about the definition of the grids, but this requires no substantial changes in the proof (the terms $2k+1$ are to be replaced by $l(x_k)$, with $l(x_{k+1})-l(x_k)=2\sgn h'(x_k)$).

\medbreak
{\bf Acknowledgements.}
I would very much like to thank Wendelin Werner for his help and advice,
especially regarding Watts' formula. I would also like to thank
Jim Pitman for stimulating discussions on intertwining relations.

-----------------------

Laboratoire de Math\'ematiques, B\^at. 425

Universit\'e Paris-Sud, F-91405 Orsay cedex, France

julien.dubedat@math.u-psud.fr

\end{document}

%% file: wedge1.pstex_t
\begin{picture}(0,0)%
\special{psfile=wedge1.pstex}%
\end{picture}%
\setlength{\unitlength}{3947sp}%
\begingroup\makeatletter\ifx\SetFigFont\undefined%
\gdef\SetFigFont#1#2#3#4#5{%
  \reset@font\fontsize{#1}{#2pt}%
  \fontfamily{#3}\fontseries{#4}\fontshape{#5}%
  \selectfont}%
\fi\endgroup%
\begin{picture}(2887,2457)(1164,-3923)
\put(4051,-2086){\makebox(0,0)[lb]{\smash{\SetFigFont{12}{14.4}{\rmdefault}{\mddefault}{\updefault}{\color[rgb]{0,0,0}$\alpha$}%
}}}
\end{picture}

%% file: wedge2.pstex_t
\begin{picture}(0,0)%
\special{psfile=wedge2.pstex}%
\end{picture}%
\setlength{\unitlength}{3947sp}%
\begingroup\makeatletter\ifx\SetFigFont\undefined%
\gdef\SetFigFont#1#2#3#4#5{%
  \reset@font\fontsize{#1}{#2pt}%
  \fontfamily{#3}\fontseries{#4}\fontshape{#5}%
  \selectfont}%
\fi\endgroup%
\begin{picture}(5412,1439)(1789,-3369)
\put(3376,-2536){\makebox(0,0)[lb]{\smash{\SetFigFont{12}{14.4}{\rmdefault}{\mddefault}{\updefault}{\color[rgb]{0,0,0}$p=\sin^2\alpha/2$}%
}}}
\put(1801,-2536){\makebox(0,0)[lb]{\smash{\SetFigFont{12}{14.4}{\rmdefault}{\mddefault}{\updefault}{\color[rgb]{0,0,0}$p$}%
}}}
\put(6526,-2611){\makebox(0,0)[lb]{\smash{\SetFigFont{12}{14.4}{\rmdefault}{\mddefault}{\updefault}{\color[rgb]{0,0,0}$p$}%
}}}
\put(7201,-2986){\makebox(0,0)[lb]{\smash{\SetFigFont{12}{14.4}{\rmdefault}{\mddefault}{\updefault}{\color[rgb]{0,0,0}$2\alpha-\pi/2$}%
}}}
\put(6526,-2086){\makebox(0,0)[lb]{\smash{\SetFigFont{12}{14.4}{\rmdefault}{\mddefault}{\updefault}{\color[rgb]{0,0,0}$p$}%
}}}
\put(5326,-2236){\makebox(0,0)[lb]{\smash{\SetFigFont{12}{14.4}{\rmdefault}{\mddefault}{\updefault}{\color[rgb]{0,0,0}$q$}%
}}}
\put(5551,-3286){\makebox(0,0)[lb]{\smash{\SetFigFont{12}{14.4}{\rmdefault}{\mddefault}{\updefault}{\color[rgb]{0,0,0}$q$}%
}}}
\put(2551,-2086){\makebox(0,0)[lb]{\smash{\SetFigFont{12}{14.4}{\rmdefault}{\mddefault}{\updefault}{\color[rgb]{0,0,0}$q=\cos^2\alpha/2$}%
}}}
\put(3076,-3211){\makebox(0,0)[lb]{\smash{\SetFigFont{12}{14.4}{\rmdefault}{\mddefault}{\updefault}{\color[rgb]{0,0,0}$q$}%
}}}
\end{picture}

%% file: wedge3.pstex_t
\begin{picture}(0,0)%
\special{psfile=wedge3.pstex}%
\end{picture}%
\setlength{\unitlength}{3947sp}%
\begingroup\makeatletter\ifx\SetFigFont\undefined%
\gdef\SetFigFont#1#2#3#4#5{%
  \reset@font\fontsize{#1}{#2pt}%
  \fontfamily{#3}\fontseries{#4}\fontshape{#5}%
  \selectfont}%
\fi\endgroup%
\begin{picture}(3387,2124)(2389,-2773)
\put(5776,-1336){\makebox(0,0)[lb]{\smash{\SetFigFont{14}{16.8}{\rmdefault}{\mddefault}{\updefault}{\color[rgb]{0,0,0}$\alpha$}%
}}}
\end{picture}

%% file: wedge4.pstex_t
\begin{picture}(0,0)%
\special{psfile=wedge4.pstex}%
\end{picture}%
\setlength{\unitlength}{3947sp}%
\begingroup\makeatletter\ifx\SetFigFont\undefined%
\gdef\SetFigFont#1#2#3#4#5{%
  \reset@font\fontsize{#1}{#2pt}%
  \fontfamily{#3}\fontseries{#4}\fontshape{#5}%
  \selectfont}%
\fi\endgroup%
\begin{picture}(5799,2015)(2314,-2514)
\put(4051,-2161){\makebox(0,0)[lb]{\smash{\SetFigFont{14}{16.8}{\rmdefault}{\mddefault}{\updefault}{\color[rgb]{0,0,0}=}%
}}}
\put(4051,-961){\makebox(0,0)[lb]{\smash{\SetFigFont{14}{16.8}{\rmdefault}{\mddefault}{\updefault}{\color[rgb]{0,0,0}=}%
}}}
\put(6301,-2161){\makebox(0,0)[lb]{\smash{\SetFigFont{14}{16.8}{\rmdefault}{\mddefault}{\updefault}{\color[rgb]{0,0,0}-}%
}}}
\put(6301,-961){\makebox(0,0)[lb]{\smash{\SetFigFont{14}{16.8}{\rmdefault}{\mddefault}{\updefault}{\color[rgb]{0,0,0}-}%
}}}
\end{picture}

%% file: wedge5.pstex_t
\begin{picture}(0,0)%
\special{psfile=wedge5.pstex}%
\end{picture}%
\setlength{\unitlength}{3947sp}%
\begingroup\makeatletter\ifx\SetFigFont\undefined%
\gdef\SetFigFont#1#2#3#4#5{%
  \reset@font\fontsize{#1}{#2pt}%
  \fontfamily{#3}\fontseries{#4}\fontshape{#5}%
  \selectfont}%
\fi\endgroup%
\begin{picture}(3975,3585)(4891,-5101)
\put(7096,-1651){\makebox(0,0)[lb]{\smash{\SetFigFont{12}{14.4}{\rmdefault}{\mddefault}{\updefault}{\color[rgb]{0,0,0}A}%
}}}
\put(4891,-4876){\makebox(0,0)[lb]{\smash{\SetFigFont{12}{14.4}{\rmdefault}{\mddefault}{\updefault}{\color[rgb]{0,0,0}B}%
}}}
\put(8866,-4861){\makebox(0,0)[lb]{\smash{\SetFigFont{12}{14.4}{\rmdefault}{\mddefault}{\updefault}{\color[rgb]{0,0,0}C}%
}}}
\put(6721,-5101){\makebox(0,0)[lb]{\smash{\SetFigFont{12}{14.4}{\rmdefault}{\mddefault}{\updefault}{\color[rgb]{0,0,0}X}%
}}}
\put(5776,-3361){\makebox(0,0)[lb]{\smash{\SetFigFont{12}{14.4}{\rmdefault}{\mddefault}{\updefault}{\color[rgb]{0,0,0}D}%
}}}
\put(8041,-3466){\makebox(0,0)[lb]{\smash{\SetFigFont{12}{14.4}{\rmdefault}{\mddefault}{\updefault}{\color[rgb]{0,0,0}E}%
}}}
\end{picture}

%% file: wedge2.bbl
\begin{thebibliography}{99}

\bibitem[AS65]{AS} M. Abramotiz, I. Stegun (eds), {\it Handbook of
    mathematical functions}, 1965
\bibitem[Ahl79]{Ahl} L. Ahlfors, {\it Complex Analysis}, 3rd edition,
  McGraw-Hill, 1979
\bibitem[Bat53]{B} H. Bateman, {\it Higher transcendental functions},
  McGraw-Hill, 1953
\bibitem[Ca92]{Ca} J.L. Cardy, {\it Critical percolation in finite
    geometries}, J. Phys. A: Math. Gen 25, pp 201--206, 1992
\bibitem[CPY98]{CPY} P. Carmona, F. Petit, M. Yor, {\it Beta-gamma
    random variables and intertwining relations between certain Markov
    processes}, Rev. Mat. Iberoamericana 14, no 2, pp 311--367, 1998 
\bibitem[Kes82]{Ke} H. Kesten, {\it Percolation theory for
    mathematicians}, Birkha\"user, 1982
\bibitem[LSW01]{LSW1}
G. Lawler, O. Schramm, W. Werner, {\it Values of Brownian intersection
  exponents. I. Half-plane exponents}, Acta Math. 187, no 2, pp
237--273, 2001
\bibitem[LSW02]{LSW} G. Lawler, O. Schramm, W. Werner, {\it Conformal
      restriction: the chordal case}, arXiv:math.PR/0209343, 2002 
\bibitem[RevYor91]{RY} D. Revuz, M. Yor, {\it Continuous Martingales and
    Brownian Motion}, Springer-Verlag, 1991 
\bibitem[RogPit81]{RP} L.C.G. Rogers, J. Pitman, {\it Markov
    functions}, Ann. Probab. 9, pp 573--582, 1981 
\bibitem[RogWil93]{RW} L.C.G. Rogers, D. Williams, {\it Diffusions,
    Markov Processes, and Martingales. Volume One: Foundations}, 2nd
  ed., John Wiley and Sons, 1993
\bibitem[Smi01]{S} S. Smirnov, {\it Critical percolation in the
    plane. I. Conformal Invariance and Cardy's formula II. Continuum
    scaling limit}, in preparation, 2001
\bibitem[VarWil85]{VW} S.R.S. Varadhan, R.J. Williams, {\it Brownian
    Motion in a Wedge with Oblique Reflection},
  Comm. Pure. App. Math., 38, pp 405--443, 1985 
\bibitem[V03]{V} B. Vir\'ag, {\it Brownian beads}, in preparation 
\bibitem [Wa96]{Wa} G.M.T. Watts, {\it A crossing probability for
    critical percolation in two dimensions}, J. Phys. A:
  Math. Gen. 29, pp 363--368, 1996
\bibitem[Wer01]{W} W. Werner, {\it Critical exponents, conformal
    invariance and planar Brownian motion}, Proceedings of the 3rd
  Europ. Congress Math., Prog. Math. 202, Birkh\"auser, pp 87--103, 2001 


\end{thebibliography}
